\documentclass[11pt]{article}
\usepackage{amssymb}
\usepackage{amsmath}

\newtheorem{theorem}{Theorem}

\newtheorem{corollary}[theorem]{Corollary}

\newtheorem{definition}[theorem]{Definition}

\newtheorem{lemma}[theorem]{Lemma}

\newtheorem{proposition}[theorem]{Proposition}
\newtheorem{remark}[theorem]{Remark}

\newenvironment{proof}[1][Proof]{\noindent\textbf{#1.} }{\ \rule{0.5em}{0.5em}}
\input{tcilatex}
\begin{document}

\title{On the Cauchy problem for integro-differential operators in H\"{o}%
lder classes and the uniqueness of the martingale problem}
\author{R. Mikulevicius and H. Pragarauskas \\
University of Southern California, Los Angeles\\
Institute of Mathematics and Informatics, Vilnius }
\maketitle

\begin{abstract}
The existence and uniqueness in H\"{o}lder spaces of solutions of the Cauchy
problem to parabolic integro-differential equation of the order $\alpha \in
(0,2)$ is investigated. The principal part of the operator has kernel $%
m(t,x,y)/|y|^{d+\alpha }$ with a bounded nondegenerate $m,$ H\"{o}lder in $x$
and measurable in $y.$ The result is applied to prove the uniqueness of the
corresponding martingale problem.

MSC classes: 45K05, 60J75, 35B65

\end{abstract}

\section{Introduction}

In this note we consider the Cauchy problem%
\begin{eqnarray}
\partial _{t}u(t,x) &=&Lu(t,x)+f(t,x),(t,x)\in H=[0,T]\times \mathbf{R}^{d},
\label{intr1} \\
u(0,x) &=&0  \notag
\end{eqnarray}%
in H\"{o}lder spaces for a class of integrodifferential operators $L=A+B$ of
the order $\alpha \in (0,2)$ whose principal part $A$ is of the form%
\begin{eqnarray}
Au(t,x) &=&A_{t}u(t,x)  \label{0} \\
&=&\int \left[ u(x+y)-u(x)-\chi _{\alpha }(y)(\nabla u(x),y)\right] m(t,x,y)%
\frac{dy}{|y|^{d+\alpha }}  \notag
\end{eqnarray}%
with $\chi _{\alpha }(y)=1_{\alpha >1}+1_{\alpha =1}1_{\left\{ |y|\leq
1\right\} }$. We notice that the operator $A$ is the generator of an $\alpha 
$-stable process. If $m=1,$then $A=c\left( -\Delta \right) ^{\alpha /2}$
(fractional Laplacian) is the generator of a spherically symmetric $\alpha $%
-stable process. The part $B$ is a perturbing, subordinated operator.

In \cite{MiP922}, the problem was considered assuming that $m$ is Holder
continuous in $x$, homogeneous of order zero and smooth in $y$ and for some $%
\eta >0$ {\ 
\begin{equation}
\int_{S^{d-1}}|(w,\xi )|^{\alpha }m(t,x,w)\mu _{d-1}(dw)\geq \eta ,\quad
(t,x)\in H,|\xi |=1,  \label{1}
\end{equation}%
where }$\mu _{d-1}$ is the Lebesgue measure on the unit sphere $S^{d-1}$ in $%
\mathbf{R}^{d}$. In \cite{AbK09}, the existence and uniqueness of a solution
to (\ref{intr1}) in H\"{o}lder spaces was proved analytically for $m$ H\"{o}%
lder continuous in $x$, smooth in $y$ and such that 
\begin{equation}
C\geq m\geq \delta >0  \label{2}
\end{equation}%
without assumption of homogeneity in $y$. \ The elliptic problem $Lu=f$ in $%
\mathbf{R}^{d}$ was considered in \cite{bas}, \cite{cafsilv}, \cite{KimDong1}%
$.$ In \cite{cafsilv}, the interior H\"{o}lder estimates (in a non-linear
case as well) were studied assuming (\ref{2}) and $m(x,y)=m(y)=m(-y)$. In 
\cite{bas}, the apriori estimates were derived in Holder classes assuming (%
\ref{2}) and Holder continuity of $m$ in $x$, except the case $\alpha =1$.
Similar results, including the case $\alpha =1$ were proved in \cite%
{KimDong1}. The equation (\ref{intr1}) with $\alpha =1$ can be regarded as a
linearization of the quasigeostrophic equation (see \cite{cav1}).

In this note, we consider he problem (\ref{intr1}), assuming that $m$ is
measurable, Holder continuous in $x$ and 
\begin{equation}
C\geq m\geq m_{0},  \label{3}
\end{equation}
where the function $m_{0}=m_{0}(t,x,y)$ is smooth and homogeneous in $y$ and
satisfies (\ref{1}). So, the density $m$ can degenerate on a substantial set.

A certain aspect of the problem is that the symbol of the operator $A,$%
\begin{equation*}
\psi (t,x,\xi )=\int \left[ e^{i(\xi ,y)}-1-\chi _{\alpha }(y)i(\xi ,y)%
\right] m(t,x,y)\frac{dy}{|y|^{d+\alpha }}
\end{equation*}%
is not smooth in $\xi $ and the standard Fourier multiplier results (for
example, used in \cite{MiP922}) do not apply in this case. Instead we use
direct analytic and probabilistic arguments. We start with equation (\ref%
{intr1}) assuming that $B=0$, the input function $f$ is smooth and the
function $m=m(t,Y)$ is smooth and homogeneous in $y,$ satisfies (\ref{1})
and does not depend on $x$. This case of equation (\ref{intr1}) was
considered in \cite{MiP922} and \cite{MiP09} and the estimates of its
solution in H\"{o}lder spaces were derived. Then we use the Ito-Wentzell
formula to pass to $m(t,y)$ which is only measurable and satisfies (\ref{3})
and obtain a solution of (\ref{intr1}) with all the estimates retained. The
case of variable coefficients is considered by using partition of unity and
deriving apriori Schauder estimates in H\"{o}lder-Zygmund spaces. Finally,
we apply the continuation by parameter method to extend solvability of an
equation with constant coefficients to that of (\ref{intr1}).

As an application, we consider the martingale problem associated to $L$.
Since the coefficients are H\"{o}lder the existence of a martingale solution
is trivial. Applying the Ito formula to the solution of (\ref{intr1}), we
prove the weak uniqueness of the solution to the martingale problem,
generalizing so the uniqueness results in \cite{AbK09} and \cite{MiP923}$.$

The note is organized as follows. In Section 2, the main theorem is stated.
In Section 3, the essential technical results are presented. The case of the
equation with constant coefficients not depending on the spacial variable is
considered in Section 4. The main theorem is proved in Section 5. In Section
6 the uniqueness of the associated martingale problem is considered.

\section{Notation and main results}

Denote $H=[0,T]\times \mathbf{R}^{d}$, $\mathbf{N}=\{0,1,2,\ldots \}$, $%
\mathbf{R}_{0}^{d}=\mathbf{R}^{d}\backslash \{0\}$. If $x,y\in \mathbf{R}%
^{d} $, we write 
\begin{equation*}
(x,y)=\sum_{i=1}^{d}x_{i}y_{i},|x|=(x,x)^{1/2}.
\end{equation*}
For a function $u=u(t,x)$ on $H$, we denote its partial derivatives by $%
\partial _{t}u=\partial u/\partial t,\partial _{i}u=\partial u/\partial
x_{i},\partial _{ij}^{2}u=\partial ^{2}u/\partial x_{i}\partial x_{j}$ and $%
D^{\gamma }u=\partial ^{|\gamma |}u/\partial x_{i}^{\gamma _{1}}\ldots
\partial x_{d}^{\gamma _{d}},$ where multiindex $\gamma =(\gamma _{1},\ldots
,\gamma _{d})\in \mathbf{N}^{d},\nabla u=(\partial _{1}u,\ldots ,\partial
_{d}u)$ denotes the gradient of $u$ with respect to $x$.

For a function $u$ on $H$ and $\beta \in (0,1]$, we write%
\begin{eqnarray*}
|u|_{0} &=&\sup_{t,x}|u(t,x)|, \\
\lbrack u]_{\beta } &=&\sup_{t,x,h\neq 0}\frac{|u(t,x+h)-u(t,x)|}{h^{\beta }}%
\text{ if }\beta \in (0,1), \\
\lbrack u]_{\beta } &=&\sup_{t,x,h\neq 0}\frac{|u(t,x+h)+u(t,x-h)-2u(t,x)|}{%
|h|}\text{ if }\beta =1.
\end{eqnarray*}%
For $\beta =[\beta ]^{-}+\left\{ \beta \right\} ^{+}>0$, where $[\beta
]^{-}\in \mathbf{N}$ and $\left\{ \beta \right\} ^{+}\in (0,1]$, we denote $%
C^{\beta }(H)$ denote the space of measurable functions $u$ on $H$ such that
the norm 
\begin{equation*}
|u|_{\beta }=\sum_{|\gamma |\leq \lbrack \beta ]^{-}}|D^{\gamma
}u|_{0}+\sup_{|\gamma |=[\beta ]^{-}}[D^{\gamma }u]_{\left\{ \beta \right\}
^{+}}.
\end{equation*}%
Accordingly, $C^{\beta }(\mathbf{R}^{d})$ denotes the corresponding space of
functions on $\mathbf{R}^{d}$.~The classes $C^{\beta }$ coincide with H\"{o}%
lder spaces if $\beta \notin \mathbf{N}$ (see 1.2.2 of \cite{Tri92}).

For $\alpha \in (0,2)$ and $u\in C^{\alpha +\beta }(H),$ we define the
fractional Laplacian%
\begin{equation}
\partial ^{\alpha }u(t,x)=\int [u(t,x+y)-u(t,x)-\left( \nabla
u(t,x),y\right) \chi _{\alpha }(y)]\frac{dy}{|y|^{d+\alpha }},  \label{fo4}
\end{equation}%
where $\chi ^{(\alpha )}(y)=\mathbf{1}_{\{|y|\leq 1\}}\mathbf{1}_{\{\alpha
=1\}}+\mathbf{1}_{\{\alpha \in (1,2)\}}$.

We denote $C_{b}^{\infty }(H)$ the space of bounded infinitely
differentiable in $x$ functions whose derivatives are bounded.

$C=C(\cdot ,\ldots ,\cdot )$ denotes constants depending only on quantities
appearing in parentheses. In a given context the same letter is (generally)
used to denote different constants depending on the same set of arguments.

Let $(U,\mathcal{U)}$ be a measurable space with a non-negative measure $\pi
(d\upsilon )$ on it$.$

Let $\alpha \in (0,2)$ and $\beta \in (0,1]$ be fixed. Let $m:H\times 
\mathbf{R}_{0}^{d}\rightarrow \lbrack 0,\infty ),b:H\rightarrow \mathbf{R}%
^{d},c:H\times U\rightarrow \mathbf{R}^{d}$ and $\rho :H\times U\rightarrow 
\mathbf{R}$ be measurable functions. We also introduce an auxiliary function 
$m_{0}:[0,T]\times \mathbf{R}_{0}^{d}\rightarrow \lbrack 0,\infty )$ and fix
positive constants $K$ and $\mu $. Throughout the paper we assume that the
function $m_{0}$ satisfies the following conditions. \medskip 

\noindent \textbf{Assumption} $\mathbf{A}_{0}.$ (i) The function $%
m_{0}=m_{0}(t,y)\geq 0$ is measurable, homogeneous in $y$ with index zero,
differentiable in $y$ up to the order $d_{0}=[\frac{d}{2}]+1$ and%
\begin{equation*}
|D_{y}^{\gamma }m_{0}^{(\alpha )}(t,y)|\leq K
\end{equation*}%
for all $t\in \lbrack 0,T]$, $y\in \mathbf{R}_{0}^{d}$ and multiindices $%
\gamma \in \mathbf{N}_{0}^{d}$ such that $|\gamma |\leq d_{0}$;

(ii) If $\alpha =1$, then for all $t\in \lbrack 0,T]$ 
\begin{equation*}
\int_{S^{d-1}}wm_{0}(t,w)\mu _{d-1}(dw)=0,
\end{equation*}%
where $S^{d-1}$ is the unit sphere in $\mathbf{R}^{d}$ and $\mu _{d-1}$ is
the Lebesgue measure on it;

(iii) For all $t\in \lbrack 0,T]$ 
\begin{equation*}
\inf_{|\xi |=1}\int_{S^{d-1}}|(w,\xi )|^{\alpha }m_{0}(t,w)\mu
_{d-1}(dw)\geq \eta >0.
\end{equation*}

\begin{remark}
\label{r10}\emph{The nondegenerateness assumption $A_{0}$ (iii) holds with
certain $\delta >0$ if, e.g.%
\begin{equation*}
\inf_{t\in \lbrack 0,T],w\in \Gamma }m_{0}^{(\alpha )}(t,w)>0
\end{equation*}%
for a measurable subset $\Gamma \subset S^{d-1}$ of positive Lebesgue
measure.}
\end{remark}

Further we will use the following assumptions.

\textbf{Assumption }$\mathbf{A.}$ (i) For all $(t,x)\in H,y\in \mathbf{R}%
_{0}^{d},$%
\begin{equation*}
|m(\cdot ,y)|_{\beta }\leq K
\end{equation*}%
and%
\begin{equation*}
m(t,x,y)\geq m_{0}(t,y),
\end{equation*}%
where the function $m_{0}$ satisfies Assumption $\mathbf{A}_{0}$;

(ii) If $\alpha =1$, then for all $(t,x)\in H$ and $r\in (0,1),$%
\begin{equation*}
\int_{r<|y|\leq 1}ym(t,x,y)\frac{dy}{|y|^{d+\alpha }}=0.
\end{equation*}

We will assume that there is a decreasing sequence of subsets $U_{n}\in 
\mathcal{U}$ such that $U=\cup _{n}U_{n}^{c}$ and the following assumptions
hold.

\textbf{Assumption B1. }(i) for all $(t,x)\in H,$%
\begin{equation*}
\int_{U_{1}}|c(t,x,\upsilon )|^{\alpha }\pi (d\upsilon
)+\int_{U_{1}^{c}}|c(t,x,\upsilon )|^{\alpha \wedge 1}\wedge 1\pi (d\upsilon
)\leq K
\end{equation*}%
(ii) for $\alpha \in (0.2)$\textbf{\ }%
\begin{equation*}
\lim_{\varepsilon \rightarrow 0}\sup_{t,x}\int 1_{|c(t,x,\upsilon )|\leq
\varepsilon }|c(t,x,\upsilon )|^{\alpha }\pi (d\upsilon )=0
\end{equation*}

\textbf{Assumption B2}. (i) 
\begin{equation*}
|b|_{\beta }+|l|_{\beta }\leq K;
\end{equation*}

(ii) If $\alpha \in \lbrack 1,2)$, then there is a constant $C$ such that
for all $(t,x)\in H,h\in \mathbf{R}^{d},$%
\begin{equation*}
\int_{U_{1}}|c(t,x,\upsilon )-c(t,x+h,\upsilon )|^{\alpha }\pi (d\upsilon
)\leq C|h|^{\alpha \beta },
\end{equation*}%
and%
\begin{equation*}
\int_{U_{1}^{c}}[|c(t,x,\upsilon )-c(t,x+h,\upsilon )|\wedge 1]\pi
(d\upsilon )\leq C|h|^{\beta };
\end{equation*}

(iii) If $\alpha <1,$ then there is $\beta ^{\prime }\,\ $such that $\alpha
+\beta >\alpha +\beta ^{\prime }\geq \beta $ and there is a constant $C$
such that for all $(t,x),\in H,h\in \mathbf{R}^{d},$ 
\begin{eqnarray*}
\int_{U_{1}}|c(t,x,\upsilon )-c(t,x+h,\upsilon )|^{(\alpha +\beta ^{\prime
})\wedge 1}\pi (d\upsilon ) &\leq &C|h|^{\beta }, \\
\int_{U_{1}^{c}}|c(t,x,\upsilon )-c(t,x+h,\upsilon )|^{(\alpha +\beta
^{\prime })\wedge 1}\wedge 1\pi (d\upsilon ) &\leq &C|h|^{\beta };
\end{eqnarray*}

(iv) For all $\upsilon \in U,$%
\begin{equation*}
|\rho \left( \cdot ,\upsilon \right) |_{\beta }\leq K.
\end{equation*}

For $(t,z)\in H,u\in C^{\alpha +\beta }(\mathbf{R}^{d})$ we introduce the
operators%
\begin{equation*}
A_{t,z}u(x)=A_{t,z}^{m}u(x)=\int_{\mathbf{R}^{d}}[u(x+y)-u(x)-(\nabla
u(x),y)\chi _{\alpha }(y)]m(t,z,y)\frac{dy}{|y|^{d+\alpha }},
\end{equation*}

\begin{eqnarray*}
B_{t,z,\bar{z}}u(x) &=&(b(t,z)\nabla u(x))1_{1\leq \alpha
<2}+\int_{U}[u(x+c(t,z,\upsilon ))-u(x) \\
&&-(\nabla u(x),c(t,z,\upsilon ))1_{U_{1}}(\upsilon )1_{1<\alpha <2}]\rho (t,%
\bar{z},\upsilon )\pi (d\upsilon ) \\
&&+l(t,z)u(x),
\end{eqnarray*}%
and

\begin{equation}
L_{t,z}u(x)=A_{t,z}u(x)+B_{t,z,z}u(x).  \label{for0}
\end{equation}%
For brevity of notation, we write%
\begin{eqnarray}
Au(t,x) &=&A_{t}u(x)=A_{t,x}u(x),Bu(t,x)=B_{t}u(x)=B_{t,x,x}u(x),
\label{for2} \\
Lu(t,x) &=&L_{t}u(x)=L_{t,x}u(x),L=A+B.  \notag
\end{eqnarray}%
According to Assumptions \textbf{A, B1, B2}, the operator $A$ represents the
principal part of $L$ and the operator $B$ is a lower order operator.

\begin{remark}
A simple example of $U,U_{n},\pi (d\upsilon )$ is $U=\mathbf{R}%
_{0}^{d},U_{n}=\{\upsilon :|\upsilon |<1/n\},c(t,x,\upsilon )=\upsilon ,\pi
(d\upsilon )=d\upsilon /|\upsilon |^{d+\alpha ^{\prime }},\alpha ^{\prime
}<\alpha ,$ and%
\begin{equation*}
Bu(t,x)=\int_{\mathbf{R}_{0}^{d}}[u(x=y)-u(x)-(\nabla u(x),y)\mathbf{1}%
_{\left\{ |y|\leq 1\right\} }\mathbf{1}_{1\leq \alpha ^{\prime }<2}]\rho
(t,x,y)\frac{dy}{|y|^{d+\alpha ^{\prime }}}.
\end{equation*}
\end{remark}

For a fixed $\alpha \in (0,2),\beta \in (0,1)$ we consider the following
Cauchy problem%
\begin{eqnarray}
\partial _{t}u(t,x) &=&(L-\lambda )u(t,x)+f(t,x),(t,x)\in H,  \label{eq1} \\
u(0,x) &=&0,x\in \mathbf{R}^{d},  \notag
\end{eqnarray}%
in Holder classes $C^{\alpha +\beta }(H)$, where $\lambda \geq 0$ and $f\in
C^{\beta }(H)$.

\begin{definition}
Let $f$ be a bounded measurable function on $H.$ We say that $u\in C^{\alpha
+\beta }(H)$ is a solution of (\ref{eq1}), if for each $(t,x)\in H,$%
\begin{equation}
u(t,x)=\int_{0}^{t}[Lu(s,x)-\lambda u(s,x)+f(s,x)]ds.  \label{defs}
\end{equation}
\end{definition}

If Assumptions \textbf{A} and \textbf{B1} are satisfied, then $Lu$ is
bounded (see Proposition \ref{prop2} and Lemma \ref{l7} below). So, (\ref%
{defs}) is well defined.

The main result of the paper is the following theorem.

\begin{theorem}
\label{main}Let $\alpha \in (0,2),\beta \in (0,1]$ and Assumptions \textbf{%
A, B1, }and \textbf{B2 }be satisfied.

Then for any $f\in C^{\beta }(H)$ there exists a unique solution $u\in
C^{\alpha +\beta }(H)$ to (\ref{eq1}). Moreover, there is a constant $%
C=C(\alpha ,\beta ,d,K,\mu )$ such that%
\begin{equation*}
|u|_{\alpha +\beta }\leqslant C|f|_{\beta },
\end{equation*}%
and for all $s\leq t\leq T,$ 
\begin{equation*}
|u(t,\cdot )-u(s,\cdot )|_{\frac{\alpha }{2}+\beta }\leq
C(t-s)^{1/2}|f|_{\beta }.
\end{equation*}
\end{theorem}

\section{Auxiliary results}

We will use the following equality for the H\"{o}lder norm estimates.

\begin{lemma}
\label{r1}$($Lemma 2.1 in \textup{\cite{Kom84}}$)$ For $\delta \in (0,1)$
and $u\in C_{0}^{\infty }(\mathbf{R}^{d})$, 
\begin{equation}
u\left( x+y\right) -u(x)=C\int k^{(\delta )}(y,z)\partial ^{\delta }u(x-z)dz,
\label{22}
\end{equation}%
where the constant $C=C(\delta ,d)$ and 
\begin{equation*}
k^{(\delta )}(y,z)=|z+y|^{-d+\delta }-|z|^{-d+\delta }.
\end{equation*}%
Moreover, there is a constant $C=C(\delta ,d)$ such that for each $y\in 
\mathbf{R}^{d}$ 
\begin{equation*}
\int |k^{(\delta )}(y,z)|dz\leq C|y|^{\delta }.
\end{equation*}
\end{lemma}

The following Lemmas \ref{le4}, \ref{l8} are deterministic counterparts of
the statements proved in \cite{MiP09}.

\begin{lemma}
\label{le4}(see Corollary 15 in \cite{MiP09})Let $\beta \in (0,1],f\in
C^{\beta }(H)$. Then there is a sequence $f_{n}\in C_{b}^{\infty }(H)$ such
that%
\begin{equation*}
|f_{n}|_{\beta }\leq 2|f|_{\beta },|f|_{\beta }\leq \lim
\inf_{n}|f_{n}|_{\beta },
\end{equation*}%
and for any $0<\beta ^{\prime }<\beta $%
\begin{equation*}
|f_{n}-f|_{\beta ^{\prime }}\rightarrow 0\text{ as }n\rightarrow \infty .
\end{equation*}
\end{lemma}

\begin{proof}
The proof in \cite{MiP09} (Corollaries 13 and 15) for $\beta \in (0,1)$
covers without any changes the case $\beta =1$ as well.
\end{proof}

\begin{lemma}
\label{rem1}$($see Theorem 6.3.2 in\cite{bergl}$)$ For $\alpha \in (0,2)$, $%
\beta >0$, the norms $|u|_{\alpha ,\beta }=|u|_{0}+|\partial ^{\alpha
}u|_{\beta }$, and $|u|_{\alpha +\beta }$ are equivalent in $C^{\alpha
+\beta }$.
\end{lemma}

Let us introduce an operator $A^{0}$ defined as operator $A$ with $m$
replaced by $m_{0}$. In terms of Fourier transforms, for $u\in C_{b}^{\infty
}(H),$%
\begin{equation*}
\mathcal{F}\left( A^{0}u\right) (t,\xi )=\psi _{0}(t,\xi )\mathcal{F}u(t,\xi
),
\end{equation*}%
where%
\begin{eqnarray*}
\psi _{0}(t,\xi ) &=&-C\int_{S^{d-1}}|(w,\xi )|^{\alpha }[1-i(\tan \frac{%
\alpha \pi }{2}\text{sgn}(w,\xi )1_{\alpha \neq 1} \\
&&-\frac{2}{\pi }\text{sgn}(w,\xi )\ln |(w,\xi )|1_{\alpha =1}]m_{0}(t,w)\mu
_{d-1}(dw)
\end{eqnarray*}%
and the constant $C=C(\alpha )>0.$ Denote%
\begin{eqnarray*}
K_{s,t}(\xi ) &=&\exp \left\{ \int_{s}^{t}\psi _{0}(r,\xi )dr\right\} ,s\leq
t, \\
G_{s,t}(x) &=&\mathcal{F}^{-1}K_{s,t},G_{s,t}^{\lambda }(x)=e^{-\lambda
(t-s)}G_{s,t}(x).
\end{eqnarray*}%
According to Assumption $\mathbf{A}_{0}$, $\int |K_{s,t}(\xi )|d\xi <\infty
,s<t$. Therefore $G_{s,t}$ is the density function of a random variable
whose characteristic function is $K_{s,t}$. Hence,%
\begin{equation}
G_{s,t}\geq 0,\int G_{s,t}(y)dy=1,s<t.  \label{4}
\end{equation}%
Let $f\in C_{b}^{\infty }(H)$ and%
\begin{equation}
R_{\lambda }f(t,x)=\int_{0}^{t}\left[ G_{s,t}^{\lambda }\ast f(s,\cdot )%
\right] (x)ds,  \label{5}
\end{equation}%
where $\ast $ denotes the convolution with respect to $x$.

\begin{lemma}
\label{l8}$($see Lemmas 7 and 17 in \cite{MiP09})Let $\alpha \in (0,2)$, $%
\beta \in (0,1]$, $f\in C_{b}^{\infty }(H)$ and Assumption $\mathbf{A}_{0}$
be satisfied. Then the Cauchy problem%
\begin{eqnarray}
\partial _{t}u(t,x) &=&A^{0}u(t,x)-\lambda u(t,x)+f(t,x),(t,x)\in H,
\label{for1} \\
u(0,x) &=&0,x\in \mathbf{R}^{d},  \notag
\end{eqnarray}%
has a unique solution $u=R_{\lambda }f\in C_{b}^{\infty }(H)$. Moreover,
there are constants $C_{1}=C(\alpha ,\ \beta ,T,\ d$,$\mu ,K)$ and $%
C_{2}=C_{2}(\alpha ,d)$ such that 
\begin{equation}
|u|_{\alpha +\beta }\leqslant C_{1}|f|_{\beta },  \label{cc1}
\end{equation}%
\begin{equation}
|u|_{\beta }\leq C_{2}(\lambda ^{-1}\wedge T)|f|_{\beta }  \label{cc2}
\end{equation}%
and for all $0\leq s\leq t\leq T$ 
\begin{equation}
|u(t,\cdot )-u(s,\cdot )|_{\alpha /2+\beta }\leq C(t-s)^{1/2}|f|_{\beta }.
\label{cc3}
\end{equation}
\end{lemma}

\begin{proof}
The statement is proved in \cite{MiP09} for $\beta \in (0,1)$ (Lemmas 7 and
17). According to Lemma 7 in \cite{MiP09}, for each $f\in C_{b}^{\infty }(H)$
there is a unique solution $u=R_{\lambda }f\in C_{b}^{\infty }(H)$.
Obviously $\partial ^{1/2}u$ solve the equation (\ref{for2}) in $%
C_{b}^{\infty }(H)$ with $\partial ^{1/2}f$ as input function. Applying the
statement with $\beta =1/2$ we have%
\begin{eqnarray*}
|\partial ^{1/2}u|_{\alpha +1/2} &\leq &C|\partial ^{1/2}f|_{1/2}, \\
|\partial ^{1/2}u|_{1/2} &\leq &C_{2}(\lambda ^{-1}\wedge T)|\partial
^{1/2}f|_{1/2}, \\
|\partial ^{1/2}u(t,\cdot )-\partial ^{1/2}u(s,\cdot )|_{\alpha /2+1/2}
&\leq &C(t-s)^{1/2}|\partial ^{1/2}f|_{1/2},s\leq t\leq T.
\end{eqnarray*}%
By Lemma \ref{rem1} (using equivalence of norms), we see that there are
constants $C_{1}=C(\alpha ,\ \beta ,T,\ d$,$\mu ,K)$ and $C_{2}=C_{2}(\alpha
,d)$ such that%
\begin{eqnarray*}
|u|_{\alpha +1} &\leq &C_{1}|f|_{1}, \\
|u|_{1} &\leq &C_{2}(\lambda ^{-1}\wedge T)|f|_{1}, \\
|u(t,\cdot )-u(s,\cdot )|_{\alpha /2+1} &\leq &C_{1}(t-s)^{1/2}|f|_{1},s\leq
t\leq T.
\end{eqnarray*}%
The statement follows immediately for $\beta =1$ by repeating the proof of
Theorem 6 in \cite{MiP09} and using Lemma \ref{le4} with $\beta =1.$
\end{proof}

Let $c_{i}:U\rightarrow \mathbf{R}^{d},i=1,2,$ be measurable functions and $%
\nu (d\upsilon )$ be a $\sigma $-finite signed measure on $(U,\mathcal{U}).$
Consider the operators ($i=1,2$)%
\begin{equation*}
L^{i}=L^{c_{i}}u(x)=\int_{U}[u(x+c_{i}(\upsilon ))-u(x)-\mathbf{1}_{\alpha
\in (1,2)}\mathbf{1}_{U_{1}}(\upsilon )(\nabla u(x),c_{i}(\upsilon ))]\nu
(d\upsilon ),
\end{equation*}%
where $U_{1}\in \mathcal{U},|\nu |(U_{1}^{c})<\infty $ ($|\nu |$ is the
total variation of $\nu $).

\begin{lemma}
\label{le9}Let $\beta \in (0,1]$. Assume%
\begin{eqnarray*}
&&\mathbf{1}_{\alpha \in (1,2)}\{\int_{U_{1}}|c_{i}(\upsilon )|^{\alpha
}|\nu |(d\upsilon )+\int_{U_{1}^{c}}|c_{i}(\upsilon )|\wedge 1|\nu
|(d\upsilon )\} \\
&&+\mathbf{1}_{\alpha \in (0,1]}\int_{U}|c_{i}(\upsilon )|^{\alpha }\wedge
1|\nu |(d\upsilon ) \\
&\leq &K_{1},i=1,2.
\end{eqnarray*}

Then there is $\beta ^{\prime }\in (0,\beta )$ and a constant $C$ such that
for each $\kappa \in (0,1)$%
\begin{equation*}
\sup_{x}|L^{1}u(x)|\leq C|u|_{\alpha +\beta ^{\prime }}K_{1},
\end{equation*}%
and%
\begin{eqnarray*}
\left[ L^{1}u\right] _{\beta } &\leq &C|u|_{a+\beta }[1_{\alpha \in
(1,2)}\int_{U_{1},|c_{1}|\leq \kappa }|c_{1}|^{\alpha }d|\nu |+\mathbf{1}%
_{\alpha \in (0,1]}\int_{|c_{1}|\leq \kappa }|c_{1}|^{\alpha }d|\nu |] \\
&&+|u|_{\alpha +\beta ^{\prime }}\kappa ^{-\alpha }K_{1}.
\end{eqnarray*}

Also, there is $\beta ^{\prime }\in (0,\beta )$ such that $\alpha +\beta
^{\prime }\geq \beta $ and%
\begin{eqnarray*}
&&\sup_{x}|L^{1}u(x)-L^{2}u(x)| \\
&\leq &C|u|_{\alpha +\beta ^{\prime }}\{\mathbf{1}_{\alpha \in (0,1]}\int
|c_{1}-c_{2}|^{(\alpha +\beta ^{\prime })\wedge 1}\wedge 1d|\nu | \\
&&+\mathbf{1}_{\alpha \in (1,2)}[\int_{U_{1}^{c}}|c_{1}-c_{2}|\wedge 1d|\nu
|+(\int_{U_{1}}|c_{1}-c_{2}|^{\alpha }d|\nu |)^{1/\alpha }]\}.
\end{eqnarray*}
\end{lemma}

\begin{proof}
If $\alpha \in (0,1]$, then for any $\beta ^{\prime }\in (0,\beta )$, 
\begin{equation*}
\sup_{x}|L^{1}u(x)|\leq C|u|_{\alpha +\beta ^{\prime }}\left( \int
|c_{1}|^{\alpha }\wedge 1d|\nu |\right) .
\end{equation*}

If $\alpha \in (1,2)$, then 
\begin{equation*}
L^{1}u=\int_{U_{1}}...+\int_{U_{1}^{c}}...=L_{1}^{1}u+L_{2}^{2},
\end{equation*}%
and%
\begin{eqnarray*}
\sup_{x}|L_{1}^{1}u(x)| &\leq &C|u|_{\alpha }\int_{U_{1}}|c_{1}|^{\alpha
}d|\nu |, \\
\sup_{x}|L_{2}^{2}u(x)| &\leq &C|u|_{\alpha }\int_{U_{1}^{c}}|c_{1}|\wedge
1d|\nu |.
\end{eqnarray*}%
If $\alpha \in (0,1],$ then for each $\kappa \in (0,1)$%
\begin{equation*}
L^{1}u=\int_{|c_{1}|\leq \kappa }...+\int_{|c_{1}|>\kappa
}=L_{1}^{1}u+L_{2}^{1}u
\end{equation*}%
and 
\begin{equation*}
|L_{1}^{1}u|_{\beta }\leq C|u|_{\alpha +\beta }\int_{|c_{1}|\leq \kappa
}|c_{1}|^{\alpha }d|\nu |
\end{equation*}%
\begin{eqnarray*}
|L_{2}^{2}u|_{\beta } &\leq &C|u|_{\beta }|\nu |\left( |c_{1}|>\kappa
\right)  \\
&\leq &C\kappa ^{-\alpha }|u|_{\beta }\int |c_{1}|^{\alpha }\wedge 1d|\nu |.
\end{eqnarray*}%
If $\alpha \in (1,2),$ then%
\begin{equation*}
L^{1}u=\int_{U_{1}}...+\int_{U_{1}^{c}}...=L_{1}^{1}u+L_{2}^{1}u
\end{equation*}%
and%
\begin{eqnarray*}
L_{1}^{1}u(x) &=&\int_{U_{1},|c|\leq \kappa }\int_{0}^{1}(\partial ^{\alpha
-1}\nabla u(x-z)k^{(\alpha -1)}(z,sc_{1}),c_{1})dzd\nu  \\
+\int_{U_{1},|c_{1}|>\kappa }... &=&L_{11}^{1}u(x)+L_{12}^{1}u(x).
\end{eqnarray*}%
We have%
\begin{eqnarray*}
\lbrack L_{11}^{1}u]_{\beta } &\leq &C|\partial ^{\alpha -1}\nabla u|_{\beta
}\int_{U_{1},|c_{1}|\leq \kappa }|c_{1}|^{\alpha }d|\nu |, \\
\lbrack L_{12}^{1}u]_{\beta } &\leq &C\kappa ^{-\alpha }|u|_{1+\beta
}\int_{U_{1}}|c_{1}|^{\alpha }d|\nu |.
\end{eqnarray*}%
Also, 
\begin{eqnarray*}
L_{2}^{1}u &=&\int_{U_{1}^{c},|c_{1}|>1}...+\int_{U_{1}^{c},|c_{1}|\leq 1}...
\\
&=&\int_{U_{1}^{c},|c_{1}|\leq 1}\int_{0}^{1}(\nabla u(x+sc_{1}),c_{1})d\nu
+\int_{U_{1}^{c},|c_{1}|>1}... \\
&=&L_{21}^{1}u+L_{22}^{1}u,
\end{eqnarray*}%
and%
\begin{eqnarray*}
\lbrack L_{21}^{1}u]_{\beta } &\leq &C|\nabla u|_{\beta }\int |c_{1}|\wedge
1d|\nu |, \\
\lbrack L_{22}^{1}u]_{\beta } &\leq &C|u|_{\beta }\int |c_{1}|\wedge 1d|\nu
|.
\end{eqnarray*}%
Finally, 
\begin{eqnarray*}
|L^{1}u(x)-L^{2}u(x)| &\leq &\mathbf{1}_{\alpha \in (0,1]}\int
|u(x+c_{1})-u(x+c_{2})|d|\nu | \\
&&+\mathbf{1}_{\alpha \in
(1,2)}[\int_{U_{1}^{c}}|u(x+c_{1})-u(x+c_{2})|d|\nu | \\
&&+\int_{U_{1}}\int_{0}^{1}|(\nabla u(x+sc_{1})-\nabla u(x),c_{1}) \\
&&-(\nabla u(x+sc_{2})-\nabla u(x),c_{2})|d|\nu |].
\end{eqnarray*}%
So, there is $\beta ^{\prime }\in (0,\beta )$ such that $\alpha +\beta
^{\prime }\geq \beta $ and%
\begin{eqnarray*}
&&|L^{1}u(x)-L^{2}u(x)| \\
&\leq &C|u|_{\alpha +\beta ^{\prime }}[\mathbf{1}_{\alpha \in (0,1]}\int
|c_{1}-c_{2}|^{(\alpha +\beta ^{\prime })\wedge 1}\wedge 1d|\nu |+\mathbf{1}%
_{\alpha \in (1,2)}\int_{U_{1}^{c}}|c_{1}-c_{2}|\wedge 1d|\nu |] \\
&&+\mathbf{1}_{\alpha \in (1,2)}|\nabla u|_{\alpha
-1}[\int_{U_{1}}|c_{1}-c_{2}|^{\alpha -1}|c_{1}|d|\nu
|+\int_{U_{1}}|c_{2}|^{\alpha -1}|c_{1}-c_{2}|d|\nu |]
\end{eqnarray*}%
and the last term can be estimated by H\"{o}lder inequality.
\end{proof}

\begin{lemma}
\label{l7}(cf. Lemma 23 in \cite{MiP09})Let $\beta \in (0,1]$ and
assumptions B1-B2 be satisfied. Then for each $\varepsilon >0$ there is a
constant $C_{\varepsilon }$ such that for any $u\in C^{\alpha +\beta }(%
\mathbf{R}^{d}),t\in \lbrack 0,T],$%
\begin{equation*}
|Bu|_{\beta }\leq \varepsilon |u|_{\alpha +\beta }+C_{\varepsilon }|u|_{0}.
\end{equation*}
\end{lemma}

\begin{proof}
Let $\beta \in (0,1)$. Since for $(t,x)\in H,z,h\in \mathbf{R}^{d},$%
\begin{eqnarray*}
B_{t,x+h,z}u(x+h)-B_{t,x,z}u(x) &=&B_{t,x+h,z}u(x+h)-B_{t,x+h,z}u(x) \\
&&+B_{t,x+h,z}u(x)-B_{t,x,z}u(x),
\end{eqnarray*}%
it follows by Lemma \ref{le9} that for each $\varepsilon >0$ there is a
constant $C_{\varepsilon }$ such that for all $z\in \mathbf{R}^{d},$%
\begin{equation}
|B_{\cdot ,z}u|_{\beta }\leq \varepsilon |u|_{\alpha +\beta }+C_{\varepsilon
}|u|_{0}.  \label{nulis}
\end{equation}%
Let $\beta =1$. Since for $(t,x)\in H,z,h\in \mathbf{R}^{d},$%
\begin{eqnarray*}
&&B_{t,x+h,z}u(x+h)-2B_{t,x,z}u(x)+B_{t,x-h,z}u(x-h) \\
&=&[B_{t,x+h,z}u(x+h)-B_{t,x,z}u(x+h)]+[B_{t,x-h,z}u(x+h)-B_{t,x,z}u(x+h)] \\
&&+[B_{t,x-h,z}u(x-h)-B_{t,x,z}u(x-h)]-[B_{t,x-h,z}u(x+h)-B_{t,x,z}u(x+h)] \\
&&+B_{t,x,z}u(x+h)-2B_{t,x,z}u(x)+B_{t,x,z}u(x-h)
\end{eqnarray*}%
it follows again by Lemma \ref{le9} that for each $\varepsilon >0$ there is
a constant $C_{\varepsilon }$ such that for all $z\in \mathbf{R}^{d},$%
\begin{equation}
|B_{\cdot ,z}u|_{1}\leq \varepsilon |u|_{\alpha +1}+C_{\varepsilon }|u|_{0}.
\label{vienas}
\end{equation}%
Finally, if $\beta \in (0,1)$, then for all $(t,x)\in H,z\in \mathbf{R}^{d},$%
\begin{eqnarray*}
&&B_{t,x+h,x+h}u(x+h)-B_{t,x,x}u(x) \\
&=&B_{t,x+h,x+h}u(x+h)-B_{t,x,x+h}u(x) \\
&&+B_{t,x,x+h}u(x)-B_{t,x,x}u(x),
\end{eqnarray*}%
and%
\begin{eqnarray*}
&&|B_{t,x+h,x+h}u(x+h)-B_{t,x,x+h}u(x)| \\
&\leq &|h|^{\beta }\sup_{z}[B_{\cdot ,z}u]_{\beta }.
\end{eqnarray*}%
So, the statement follows by (\ref{nulis}) and Lemma \ref{le9}.

If $\beta =1$, then $(t,x)\in H,z,h\in \mathbf{R}^{d},$%
\begin{eqnarray*}
&&B_{t,x+h,x+h}u(x+h)-2B_{t,x,x}u(x)+B_{t,x-h,x-h}u(x-h) \\
&=&\{B_{t,x+h,x+h}u(x+h)-2B_{t,x+h,x}u(x+h)+B_{t,x+h,x-h}u(x+h)\} \\
&&+%
\{[B_{t,x-h,x-h}u(x-h)-B_{t,x-h,x}u(x-h)]-[B_{t,x+h,x-h}u(x+h)-B_{t,x+h,x}u(x+h)]\}
\\
&&+\{B_{t,x+h,x}u(x+h)-2B_{t,x,x}u(x)+B_{t,x-h,x}u(x-h)\}
\end{eqnarray*}%
and the statement follows by (\ref{vienas}) and Lemma \ref{le9}.
\end{proof}

Let $n:\mathbf{R}_{0}^{d}\rightarrow \mathbf{R}$ be a measurable function
satisfying the following conditions:

(i) there is a constant $k_{1}$ such that for all $(t,x)\in H,y\in \mathbf{R}%
_{0}^{d},$%
\begin{equation}
|n(y)|\leq k_{1};  \label{6}
\end{equation}

(ii) if $\alpha =1$, then for all $r\in (0,1),$%
\begin{equation}
\int_{r<|y|\leq 1}yn(y)\frac{dy}{|y|^{d+1}}=0.  \label{63}
\end{equation}%
For $u\in C^{\alpha +\beta }(\mathbf{R}^{d})$, we introduce the operators%
\begin{equation*}
\mathcal{A}u(x)=\int_{\mathbf{R}^{d}}\nabla _{y}^{\alpha }u(x)n(y)\frac{dy}{%
|y|^{d+\alpha }}.
\end{equation*}

\begin{proposition}
\label{prop2}Let $\alpha \in (0,2)$, $\beta \in (0,1],\beta ^{\prime }\in
(0,\beta )$ and (\ref{6}), (\ref{63}) be satisfied.

Then there are constants $C_{1}=C_{1}(\alpha ,\beta ^{\prime
},d),C_{2}=C_{2}(\alpha ,\beta ,d)$ such that for all $u\in C^{\alpha +\beta
}(H)$%
\begin{eqnarray*}
\sup_{x}|\mathcal{A}u(x)| &\leq &C_{1}k_{1}|u|_{\alpha +\beta ^{\prime }}, \\
\lbrack \mathcal{A}u]_{\beta } &\leq &C_{2}k_{1}|u|_{\alpha +\beta }.
\end{eqnarray*}
\end{proposition}

\begin{proof}
For $\alpha \in (0,1)$, let $\beta ^{\prime }\in (0,\beta )$ be such that $%
\alpha +\beta ^{\prime }<1$. Then for $u\in C^{\alpha +\beta }(\mathbf{R}%
^{d}),$ 
\begin{equation*}
\sup_{x}|\mathcal{A}u(x)|\leq Ck_{1}|u|_{a+\beta ^{\prime }}\int
(|y|^{\alpha +\beta ^{\prime }}\wedge 1)\frac{dy}{|y|^{d+\alpha }}\leq
CK|u|_{\alpha +\beta ^{\prime }}.
\end{equation*}%
For $\alpha \in \lbrack 1,2),$%
\begin{eqnarray*}
\mathcal{A}u(x) &=&\int_{|y|\leq 1}\left( \int_{0}^{1}(\nabla u(x+sy)-\nabla
u(x),y\right) ds]m(y)\frac{dy}{|y|^{d+\alpha }} \\
&&+\int_{|y|>1}[u(x+y)-u(x)-\mathbf{1}_{\alpha \in (1,2)}\left( \nabla
u(x),y\right) ]m(y)\frac{dy}{|y|^{d+\alpha }} \\
&=&L_{1}u(x)+L_{2}u(x).
\end{eqnarray*}

Obviously, for any $\mu \in (0,1)$ such that $1+\mu \in (\alpha ,\alpha
+\beta )$%
\begin{equation*}
\sup_{x}|L_{1}u(x)|\leq CK|\nabla u|_{\mu }\int_{|y|\leq 1}|y|^{1+\mu
-d-\alpha }dy\leq CK|u|_{1+\mu },
\end{equation*}%
and%
\begin{equation*}
\sup_{x}|L_{2}u(x)|\leq CK\sup_{x}(|u(x)|+|\nabla u(x)|).
\end{equation*}

In order to estimate the differences, first we note that for $u\in C^{1+\mu
}(\mathbf{R}^{d}),\mu \in (0,1),$%
\begin{eqnarray*}
&&u(x+h)-u(x)+u(x-h)-u(x) \\
&=&\int_{0}^{1}\left( \nabla u(x+sh)-\nabla u(x-sh),t\right) ds
\end{eqnarray*}%
and 
\begin{equation}
|u(x+h)-u(x)+u(x-h)-u(x)|\leq C|\nabla u|_{\mu }|h|^{1+\mu },x,h\in \mathbf{R%
}^{d}.  \label{fo3}
\end{equation}%
Also, for $u\in C^{\mu }(\mathbf{R}^{d}),\mu \in (0,1),$%
\begin{equation}
|u(x+h)-u(x)|\leq C|u|_{\mu }|h|^{\mu },x,h\in \mathbf{R}^{d}.  \label{fo5}
\end{equation}

Fix $h\in \mathbf{R}^{d}$ with $a=|h|\in (0,1).$ Then 
\begin{equation*}
\mathcal{A}u(x)=\int_{|y|\leq a}...+\int_{|y|>a}...=I_{1}(x)+I_{2}(x),x\in 
\mathbf{R}^{d},
\end{equation*}%
where%
\begin{eqnarray*}
I_{1}(x) &=&\int_{|y|\leq a}[u(x+y)-u(x)-\mathbf{1}_{\alpha \in \lbrack
1,2)}\left( \nabla u(x),y\right) ]m(y)\frac{dy}{|y|^{d+\alpha }}, \\
I_{2}(x) &=&\int_{|y|>a}[u(x+y)-u(x)-\mathbf{1}_{\alpha \in (1,2)}\left(
\nabla u(x),y\right) ]m(y)\frac{dy}{|y|^{d+\alpha }}.
\end{eqnarray*}%
For $\alpha \in (0,1),\beta \in (0,1]$, let $\beta ^{\prime }\in (0,\beta )$
and $\alpha +\beta ^{\prime }<1$. Then for $u\in C^{\alpha +\beta }(\mathbf{R%
}^{d})$ by Lemma \ref{r1},%
\begin{eqnarray*}
|I_{1}(x+h)-I_{1}(x)| &\leq &k_{1}\int_{|y|\leq a}|\partial ^{\alpha +\beta
^{\prime }}u(x+h-z)-\partial ^{\alpha +\beta ^{\prime }}u(x-z)|~|k^{(\alpha
+\beta ^{\prime })}(z,y)|\frac{dy}{|y|^{d+\alpha }} \\
&\leq &Ck_{1}|\partial ^{\alpha +\beta ^{\prime }}|_{\beta -\beta ^{\prime
}}a^{\beta -\beta ^{\prime }}\int_{|y|\leq a}|y|^{\alpha +\beta ^{\prime }}%
\frac{dy}{|y|^{d+\alpha }}\leq Ck_{1}|u|_{\alpha +\beta }a^{\beta }.
\end{eqnarray*}%
For $\alpha \in \lbrack 1,2),\beta \in (0,1]$, let $\beta ^{\prime }\in
(0,1) $ and $\alpha <1+\beta ^{\prime }<\alpha +\beta $. Then for $u\in
C^{\alpha +\beta }(\mathbf{R}^{d})$ by Lemma \ref{r1},%
\begin{eqnarray*}
&&|I_{1}(x+h)-I_{1}(x)| \\
&\leq &Ck_{1}\int_{|y|\leq a}\int_{0}^{1}\int |\left( \partial ^{\beta
^{\prime }}\nabla u(x+h-z)-\partial ^{\beta ^{\prime }}\nabla
u(x-z),y\right) |~|k^{(\beta ^{\prime })}(z,sy)|dzds\frac{dy}{|y|^{d+\alpha }%
} \\
&\leq &Ck_{1}|\partial ^{\beta ^{\prime }}\nabla u|_{\alpha +\beta -\beta
^{\prime }-1}a^{\alpha +\beta -1-\beta ^{\prime }}\int_{|y|\leq
a}|y|^{1+\beta ^{\prime }}\frac{dy}{|y|^{d+\alpha }}\leq Ck_{1}|u|_{\alpha
+\beta }a^{\beta }.
\end{eqnarray*}%
Let $\alpha \in (0,1],\beta \in (0,1),$ and $\alpha ^{\prime }<\alpha $ be
such that $\beta +\alpha -\alpha ^{\prime }<1$. By Lemma \ref{r1} and (\ref%
{fo5}),%
\begin{eqnarray*}
|I_{2}(x+h)-I_{2}(x)| &\leq &k_{1}\int_{|y|>a}|\partial ^{\alpha ^{\prime
}}u(x+h-z)-\partial ^{\alpha ^{\prime }}u(x-z)|~|k^{(\alpha ^{\prime
})}(y,z)|\frac{dy}{|y|^{d+\alpha }} \\
&\leq &Ck_{1}|\partial ^{\alpha ^{\prime }}u|_{\alpha +\beta -\alpha
^{\prime }}a^{\alpha +\beta -\alpha ^{\prime }}\int_{|y|>a}|y|^{\alpha
^{\prime }}\frac{dy}{|y|^{d+\alpha }}\leq Ck_{1}|u|_{\alpha +\beta }a^{\beta
}.
\end{eqnarray*}%
For $\alpha \in (0,1],$ $\beta =1,$ let $\alpha ^{\prime }\in (0,\alpha )$.
Then $1<\beta +\alpha -\alpha ^{\prime }<2$, 
\begin{equation*}
I_{2}(x)=\int_{|y|>a}\int \partial ^{\alpha ^{\prime }}u(x-z)k^{(\alpha
^{\prime })}(z,y)m(y)\frac{dy}{|y|^{d+\alpha }}.
\end{equation*}%
and by (\ref{fo3}),%
\begin{eqnarray*}
|I_{2}(x+h)-2I_{2}(x)+I_{2}(x-h)| &\leq &Ck_{1}|\partial ^{\alpha ^{\prime
}}u|_{\alpha +\beta -\alpha ^{\prime }}a^{\alpha +\beta -\alpha ^{\prime
}}\int_{|y|>a}|y|^{\alpha ^{\prime }}\frac{dy}{|y|^{d+\alpha }} \\
&\leq &Ck_{1}|u|_{\alpha +\beta }a^{\beta }.
\end{eqnarray*}%
Let $x,\bar{x}\in \mathbf{R}^{d},a=|x-\bar{x}|$, and $\beta ^{\prime }<1$ be
such that $\alpha +\beta ^{\prime }<2$ and $0\leq \beta -\beta ^{\prime }<1$%
. By Lemma \ref{r1}, 
\begin{eqnarray*}
&&|L_{1}u(x)-L_{1}(\bar{x})| \\
&\leq &K\int_{|y|\leq a}\int_{0}^{1}\int |\partial ^{\alpha +\beta ^{\prime
}-1}\nabla u(x-z)-\partial ^{\alpha +\beta ^{\prime }-1}\nabla u(x-z)|\times
\\
&&\times ~|k^{(\alpha +\beta ^{\prime }-1)}(sy,z)||y|\frac{dsdzdy}{%
|y|^{d+\alpha }} \\
&\leq &CK|u|_{\alpha +\beta }a^{\beta -\beta ^{\prime }}\int_{|y|\leq
a}|y|^{\alpha +\beta ^{\prime }}\frac{dy}{|y|^{d+\alpha }}=CK|u|_{\alpha
+\beta }a^{\beta }.
\end{eqnarray*}%
For $\alpha \in (1,2),\beta \in (0,1)$, let $1<\alpha ^{\prime }<\alpha $ be
such that $\alpha -\alpha ^{\prime }+\beta <1$. By Lemma \ref{r1}, 
\begin{eqnarray*}
&&|I_{2}(x+h)-I_{2}(x)| \\
&=&|\int_{|y|>a}\int_{0}^{1}\int \left( \partial ^{\alpha ^{\prime
}-1}\nabla u(x+h-z)-\partial ^{\alpha ^{\prime }-1}\nabla u(x-z),y\right)
k^{(\alpha ^{\prime }-1)}(z,sy)m(y)\frac{dzdsdy}{|y|^{d+\alpha }}| \\
&\leq &Ck_{1}|\partial ^{\alpha ^{\prime }-1}\nabla u|_{\alpha +\beta
-\alpha ^{\prime }}a^{\alpha +\beta -\alpha ^{\prime
}}\int_{|y|>a}|y|^{\alpha ^{\prime }}\frac{dzdsdy}{|y|^{d+\alpha }}\leq
Ck_{1}|u|_{\alpha +\beta }a^{\beta }.
\end{eqnarray*}%
For $\alpha \in (1,2),\beta =1,$ we have%
\begin{equation*}
I_{2}(x)=\int_{|y|>a}\int_{0}^{1}\left( \nabla u(x+sy)-\nabla u(x),y\right)
m(y)ds\frac{dy}{|y|^{d+\alpha }},
\end{equation*}%
and, by (\ref{fo3}),%
\begin{eqnarray*}
|I_{2}(x+h)-2I_{2}(x)+I_{2}(x-h)| &\leq &Ck_{1}|\nabla u|_{\alpha +\beta
-1}a^{\alpha }\int_{|y|>a}|y|\frac{dy}{|y|^{d+\alpha }} \\
&=&Ck_{1}|u|_{a+\beta }a.
\end{eqnarray*}

The statement follows.
\end{proof}

We will need a generalization of this statement. Let $n:\mathbf{R}^{d}%
\mathbf{\times R}_{0}^{d}\rightarrow \mathbf{R}$ be a measurable function
satisfying the following conditions:

(i) there is a constant $k_{1}$ such that for all $x\in \mathbf{R}^{d},y\in 
\mathbf{R}_{0}^{d}$%
\begin{equation}
|n(x,y)|\leq k_{1};  \label{91}
\end{equation}

(ii) for $\beta \in (0,1]$ there is a constant $k_{2}$ such that for all $%
y\in \mathbf{R}_{0}^{d}$%
\begin{equation}
\lbrack n(\cdot ,y)]_{\beta }\leq k_{2};  \label{9}
\end{equation}

(iii) if $\alpha =1$, then for all $x\in \mathbf{R}^{d},r\in (0,1),$%
\begin{equation}
\int_{r<|y|\leq 1}yn(x,y)\frac{dy}{|y|^{d+1}}=0.  \label{92}
\end{equation}%
For $u\in C^{\alpha +\beta }(\mathbf{R}^{d})$, we introduce an operator%
\begin{equation*}
\mathcal{A}u(x)=\int_{\mathbf{R}^{d}}\nabla _{y}^{\alpha }u(x)n(x,y)\frac{dy%
}{|y|^{d+\alpha }}.
\end{equation*}

\begin{corollary}
\label{cor1}Let $\alpha \in (0,2)$, $\beta \in (0,1],\beta ^{\prime }\in
(0,\beta )$ and (\ref{91})-(\ref{92}) be satisfied. Then there is a constant 
$C=C(\alpha ,\beta ,\beta ^{\prime },d)$ such that for all $u\in C^{\alpha
+\beta }(H)$%
\begin{equation*}
|\mathcal{A}u|_{\beta }\leq C[k_{1}|u|_{\alpha +\beta }+k_{2}|u|_{\alpha
+\beta ^{\prime }}].
\end{equation*}
\end{corollary}

\begin{proof}
For $u\in C^{\alpha +\beta }(\mathbf{R}^{d})$ consider%
\begin{equation*}
\mathcal{A}_{z}u(x)=\int_{\mathbf{R}^{d}}\nabla _{y}^{\alpha }u(x)n(z,y)%
\frac{dy}{|y|^{d+\alpha }},x,z\in \mathbf{R}^{d}.
\end{equation*}%
Let $\beta \in (0,1)$. Since for $h,x\in \mathbf{R}^{d},$%
\begin{eqnarray*}
\mathcal{A}_{x+h}u(x+h)-\mathcal{A}_{x}u(x) &=&\left( \mathcal{A}%
_{x+h}u(x+h)-\mathcal{A}_{x+h}u(x)\right) \\
&&+\left( \mathcal{A}_{x+h}u(x)-\mathcal{A}_{x}u(x)\right) ,
\end{eqnarray*}%
the statement follows by Proposition \ref{prop2}.

Let $\beta =1$. Since, similarly, for $h,x\in \mathbf{R}^{d},$%
\begin{eqnarray*}
&&\mathcal{A}u(x+h)-2\mathcal{A}u(x)+\mathcal{A}u(x-h) \\
&=&\{\mathcal{A}_{x-h}u(x+h)-2\mathcal{A}_{x}u(x+h)+\mathcal{A}_{x+h}u(x+h)\}
\\
&&+\{[\mathcal{A}_{x-h}u(x+t)-\mathcal{A}_{x}u(x+t)]-[\mathcal{A}%
_{x-h}u(x-t)-\mathcal{A}_{x}u(x-t)]\} \\
&&+\{\mathcal{A}_{x}u(x+h)-2\mathcal{A}_{x}u(x)+\mathcal{A}_{x}u(x-h)\},
\end{eqnarray*}%
the statement follows by Proposition \ref{prop2}.
\end{proof}

\section{Equation with coefficients independent of spatial variable}

In this section, we consider the Cauchy problem 
\begin{equation}
\left\{ 
\begin{array}{ll}
\partial _{t}u(t,x)=A_{t}u(t,x)-\lambda u(t,x)+f(t,x), & (t,x)\in H \\ 
u(0,x)=0, & x\in \mathbf{R}^{d},%
\end{array}%
\right.  \label{one}
\end{equation}%
assuming that the function $m(t,x,y)$ does not depend on $x$.

\begin{theorem}
\label{t1}Let $\alpha \in (0,2),\beta \in (0,1],m(t,x,y)=m(t,y)$ and
Assumption A be satisfied.

Then for each $f\in C^{\beta }(H)$ there is a unique solution $u\in
C^{\alpha +\beta }(H)$ to (\ref{one}). Moreover, the solution satisfies (\ref%
{cc1})-(\ref{cc3}).
\end{theorem}

\begin{proof}
\emph{Uniqueness. }Let $u^{1},u^{2}\in C^{\alpha +\beta }$ be two solutions
to (\ref{one}). Then the function $u=u^{1}-u^{2}$ satisfies (\ref{one}) with 
$f=0$. 

Let a nonnegative $\zeta \in C_{0}^{\infty }(\mathbf{R}^{d})$ be such that $%
\int \zeta dx=1$. Denote 
\begin{equation*}
\zeta _{\varepsilon }(x)=\varepsilon ^{-d}\zeta (x/\varepsilon ),x\in 
\mathbf{R}^{d},\varepsilon \in (0,1),
\end{equation*}%
and 
\begin{equation*}
u_{\varepsilon }(t,x)=u(t,\cdot )\ast \zeta _{\varepsilon }(x),(t,x)\in H%
\text{.}
\end{equation*}%
Then $u_{\varepsilon }$ solves (\ref{one}) with $f=0$. 

Let $\left( \Omega ,\mathcal{F},\mathbf{P}\right) $ be a complete
probability space with a filtration of $\sigma $-algebras $\mathbb{F}=\left( 
\mathcal{F}_{t}\right) _{t\geq 0}$ satisfying the usual conditions. We fix $%
t_{0}\in (0,T)$ and introduce an $\mathbb{F}$-adapted Poisson point measure $%
p(dt,dy)$ on $[0,t_{0}]\times \mathbf{R}_{0}^{d}$ with a compensator $%
m(t_{0}-t,y)dtdy/|y|^{d+\alpha }$. Let%
\begin{equation*}
q(dt,dy)=p(dt,dy)-m(t_{0}-t,y)\frac{dtdy}{|y|^{d+\alpha }}
\end{equation*}%
be the corresponding martingale measure and 
\begin{equation*}
X_{t}=\int_{s_{0}}^{t}\int \chi _{\alpha }(y)yq(ds,sy)+\int_{0}^{t}\int
(1-\chi _{\alpha }(y))yp(ds,dy)
\end{equation*}%
for $0\leq t\leq t_{0}.$ By Ito's formula%
\begin{eqnarray*}
u_{\varepsilon }(t_{0},x) &=&u_{\varepsilon }(t_{0},x)-\mathbf{E}%
u_{\varepsilon }(0,x+X_{t_{0}})e^{-\lambda t_{0}} \\
&=&\mathbf{E}\int_{0}^{t_{0}}e^{-\lambda t}\left[ \frac{\partial
u_{\varepsilon }}{\partial t}-Au_{\varepsilon }+\lambda u_{\varepsilon }%
\right] (t-t_{0},x+X_{t})dt=0.
\end{eqnarray*}%
Since $\varepsilon ,t_{0}$ and $x$ are arbitrary, we have $u=0.$

\emph{Existence. }First we prove the existence of a solution to (\ref{one})
for a smooth input function $f$. 

We introduce an $\mathbb{F}$-adapted Poisson measure $\bar{p}(dt,dz)$ on $%
[0,\infty )\times \mathbf{R}_{0}$ with a compensator $dtdz/z^{2}$. Let%
\begin{equation*}
\bar{q}(dt,dz)=\bar{p}(dt,dz)-\frac{dtdz}{z^{2}}
\end{equation*}%
be the corresponding martingale measure. According to Lemma 14.50 in  \cite%
{Jac79}, there is a measurable function $\bar{c}:[0,T]\times \mathbf{R}%
_{0}\rightarrow \mathbf{R}^{d}$ such that for every Borel $\Gamma \subseteq 
\mathbf{R}_{0}^{d}$%
\begin{equation*}
\int_{\Gamma }\left( m(t,y)-m_{0}(t,y)\right) \frac{dy}{|y|^{d+\alpha }}%
=\int 1_{\Gamma }(\bar{c}(t,z))\frac{dz}{z^{2}}.
\end{equation*}%
Let 
\begin{equation*}
Y_{t}=\int_{0}^{t}\int (1-\chi _{\alpha }(\bar{c}(t,z)))\bar{c}(t,z)\bar{p}%
(dz,dz)+\int_{0}^{t}\int \chi _{\alpha }(\bar{c}(t,z))c(t,z)\bar{q}(ds,dz).
\end{equation*}%
For $f\in C_{b}^{\infty }(H),$ we consider the equation%
\begin{eqnarray}
\partial _{t}u(t,x) &=&A_{\alpha }^{0}u(t,x)+f(t,x-Y_{t}),(t,x)\in H,
\label{d1} \\
u(0,x) &=&0,x\in \mathbf{R}^{d}.  \notag
\end{eqnarray}
By Lemma \ref{l8}, there is a unique solution $u\in C_{b}^{\infty }(H)$ to (%
\ref{d1}). Moreover, the solution satisfies (\ref{cc1})-(\ref{cc3}) $\mathbf{%
P}$-a.s. In addition, the solution is $\mathbb{F}$-adapted because%
\begin{equation*}
u(t,x)=\int_{0}^{t}\int G_{s,t}^{\lambda }(y)f(s,x-y-Y_{s})dyds.
\end{equation*}

Using (\ref{4}), we have for any multiindex $\gamma ,$%
\begin{equation*}
D^{\gamma }u(t,x)=\int_{0}^{t}\int G_{s,t}^{\lambda }(y)D_{x}^{\gamma
}f(s,x-y-Y_{s})dyds
\end{equation*}%
and%
\begin{equation}
\text{essup}_{\omega \in \Omega }\sup_{t,x}|D_{x}^{\gamma }u(t,x)|<\infty .
\label{a1}
\end{equation}

Let $\bar{A}$ be the operator defined as the operator $A$ with $m$ replaced
by $m-m_{0}$. According to (\ref{d1}) and the Ito-Wentzell formula (see \cite%
{mik}),%
\begin{eqnarray}
u(t,x+Y_{t})-u(0,x) &=&\int_{0}^{t}[\partial _{s}u(s,x+Y_{s})+\bar{A}%
u(s,x+Y_{s})+M_{t}  \label{a2} \\
&=&\int_{0}^{t}[Au(s,x+Y_{s})-\lambda u(s,x+Y_{s})+f(s,x)]ds+M_{t},  \notag
\end{eqnarray}%
where%
\begin{equation*}
M_{t}=\int_{0}^{t}\int [u(s,x+Y_{s-}+\bar{c}(t,z))-u(s,x+Y_{s-})]\bar{q}%
(ds,dz).
\end{equation*}

Taking expectation on both sides of (\ref{a2}) and using (\ref{a1}), we
conclude that the function $v(t,x)=\mathbf{E}u(t,x+Y_{t})$ belongs to $%
C_{b}^{\infty }(H)$ and solves (\ref{one}). Moreover, $v$ satisfies (\ref%
{cc1})-(\ref{cc3}) because $\ u$ satisfies (\ref{cc1})-(\ref{cc3}) $\mathbf{P%
}$-a.s.

Next we prove the existence of a solution to (\ref{one}) for $f\in C^{\beta
}(H)$. By Lemma \ref{le4}, there is a sequence $f_{n}\in C_{b}^{\infty }(H)$
such that%
\begin{equation}
|f_{n}|_{\beta }\leq 2|f|_{\beta },|f|_{\beta }\leq \lim
\inf_{n}|f_{n}|_{\beta },  \label{a3}
\end{equation}%
and for every $\kappa \in (0,\beta )$ 
\begin{equation}
\lim_{n\rightarrow \infty }|f_{n}-f|_{\kappa }=0.  \label{a4}
\end{equation}
According to the first part of the proof and (\ref{a3}), for each $n$ there
is a unique solution $u_{n}\in C^{\alpha +\beta }(H)$ to (\ref{one}) with $f$
replaced by $f_{n}.$ Moreover, there are constants $C_{1}=C_{1}(\alpha
,\beta ,\mu ,K,T)$ and $C_{2}=C_{2}(\alpha ,d)$ such that%
\begin{equation}
|u_{n}|_{\alpha +\beta }\leqslant C_{1}|f|_{\beta },  \label{ine0}
\end{equation}%
\begin{equation}
|u_{n}|_{\beta }\leq C_{2}(\alpha ,d)(\lambda ^{-1}\wedge T)|f|_{0,\beta ;p},
\label{a5}
\end{equation}%
and for all $s\leq t\leq T,$ 
\begin{equation}
|u_{n}(t,\cdot )-u_{n}(s,\cdot )|_{\alpha /2+\beta }\leq
C(t-s)^{1/2}|f|_{\beta }.  \label{ine2}
\end{equation}%
Fix an arbitrary $\kappa \in (0,\beta ).$ Again, by the first part of he
proof, there is a constant $C$ not depending on $n$ such that%
\begin{equation*}
|u_{n}|_{\alpha +\kappa }\leqslant C|f_{n}|_{\kappa }.
\end{equation*}%
Moreover, by Lemma \ref{le4} and (\ref{a4}),%
\begin{equation*}
|u_{n}-u_{k}|_{\alpha +\kappa }\leqslant C|f_{n}-f_{k}|_{\kappa }\rightarrow
0
\end{equation*}%
as $n,k\rightarrow \infty $. Hence, there is $u\in C^{\alpha +\kappa }(H)$
such that%
\begin{equation}
\lim_{n\rightarrow \infty }|u_{n}-u|_{\alpha +\kappa }=0.  \label{a6}
\end{equation}%
According to (\ref{ine0}) and (\ref{a6}), we have $[\partial ^{\alpha
}u_{n}]_{\beta }\leq C_{1}|f|_{\beta }$ and \TEXTsymbol{\vert}$\partial
^{\alpha }u_{n}-\partial ^{\alpha }u|_{0}\rightarrow 0$ as $n\rightarrow
\infty $. Therefore $\left[ \partial ^{\alpha }u\right] _{\beta }\leq
C_{1}|f|_{\beta }$ and $u\in C^{\alpha +\beta }(H)$. Passing to the limit in
(\ref{ine0})-(\ref{ine2}) as $n\rightarrow \infty $ we conclude that $u$
satisfies (\ref{cc1})-(\ref{cc3}). Finally, passing to the limit in the
equation%
\begin{equation*}
u_{n}(t,x)=\int_{0}^{t}\left[ Au_{n}-\lambda u_{n}+f_{n}\right] (s,x)ds
\end{equation*}%
and using Corollary \ref{cor1}, we conclude that $u$ solves (\ref{one}). 

The theorem is proved.
\end{proof}

\section{Proof of Theorem \protect\ref{main}}

We follow the proof of Theorem 5 in \cite{MiP09} with obvious changes.

It is well known that for an arbitrary but fixed $\delta >0$ there is a
family of cubes $D_{k}\subseteq \tilde{D}_{k}\subseteq \mathbf{R}^{d}$ and a
family of deterministic functions $\eta _{k}\in C_{0}^{\infty }(\mathbf{R}%
^{d})$ with the following properties:

1. For all $k\geq 1,D_{k}$ and $\tilde{D}_{k}$ have a common center $x_{k},$
diam $D_{k}\leq \delta ,$dist$(D_{k},\mathbf{R}^{d}\backslash \tilde{D}%
_{k})\leq C\delta $ for a constant $C=C(d)>0,\cup _{k}D_{k}=\mathbf{R}^{d},$
and $1\leq \sum_{k}1_{\tilde{D}_{k}}\leq 2^{d}.$

2. For all $k$, $0\leq \eta _{k}\leq 1,\eta _{k}=1$ in $D_{k},\eta _{k}=0$
outside of $\tilde{D}_{k}$ and for all multiindices $\gamma ,|\gamma |\leq 3,
$%
\begin{equation*}
|\partial ^{\gamma }\eta _{k}|\leq C(d)\delta ^{-|\gamma |}.
\end{equation*}%
For $\alpha \in (0,2),k\geq 1,$ denote%
\begin{eqnarray*}
A_{k}u(t,x) &=&A_{t,x_{k}}u(t,x), \\
E_{k}u(t,x) &=&\int [u(t,x+y)-u(t,x)][\eta _{k}(x+y)-\eta
_{k}(x)]m(t,x_{k},y)\frac{dy}{|y|^{d+\alpha }}, \\
E_{k,1}u(t,x) &=&\int [u(t,x+y)-u(t,x)][\eta _{k}(x+y)-\eta _{k}(x)]\frac{dy%
}{|y|^{d+\alpha }}, \\
F_{k}u(t,x) &=&u(t,x)A_{k}\eta _{k}(x),F_{k,1}u(t,x)=u(t,x)\partial ^{\alpha
}\eta _{k}(x).
\end{eqnarray*}%
We notice that 
\begin{equation}
A_{k}(u\eta _{k})=\eta _{k}A_{k}u+E_{k}u+F_{k}u  \label{f0}
\end{equation}%
and%
\begin{equation}
\partial ^{\alpha }(u\eta _{k})=\eta _{k}\partial ^{\alpha
}u+E_{k,1}u+F_{k,1}u.  \label{f00}
\end{equation}

It is readily checked that there is $\beta ^{\prime }\in (0,\beta )$ and a
constant $C=$ $C=C(\alpha ,\beta ,d,K,\delta )$ such that%
\begin{equation*}
\sup_{k}\left( |E_{k}^{(\alpha )}u(t,\cdot )|_{\beta }+|E_{k,1}^{(\alpha
)}u(t,\cdot )|_{\beta }\right) \leq C|u|_{\alpha +\beta ^{\prime }}
\end{equation*}%
and, by Corollary \ref{cor1},%
\begin{equation}
\sup_{k}\left( |F_{k}^{(\alpha )}u(t,\cdot )|_{\beta }+|F_{k,1}^{(\alpha
)}u(t,\cdot )|_{\beta }\right) \leq C|u|_{\beta }.  \label{f2}
\end{equation}%
Hence, for each $\varepsilon >0$ there exists a constant $C=C(\alpha ,\beta
,d,K,\delta ,\varepsilon )$ such that%
\begin{equation}
\sup_{k}\left( |E_{k}^{(\alpha )}u(t,\cdot )|_{\beta }+|E_{k,1}^{(\alpha
)}u(t,\cdot )|_{\beta }\right) \leq \varepsilon |u|_{\alpha +\beta
}+C|u|_{0}.  \label{f20}
\end{equation}

Elementary calculation shows that for every $u\in C^{\alpha +\beta }(H),$%
\begin{eqnarray}
|u|_{0} &\leq &\sup_{k}\sup_{x}|\eta _{k}(x)u(x)|,  \notag \\
|u|_{\beta } &\leq &\sup_{k}|\eta _{k}u|_{\beta }+C|u|_{0},  \label{f3} \\
\sup_{k}|\eta _{k}u|_{\beta } &\leq &|u|_{\beta }+C|u|_{0},  \notag
\end{eqnarray}
the constant $C=C(\beta ,d,\delta ).$ By (\ref{f00}) and (\ref{f3}), we have 
\begin{eqnarray*}
|u|_{\alpha ,\beta } &=&|u|_{0}+|\partial ^{\alpha }u|_{\beta }\leq
\sup_{k}|\eta _{k}\partial ^{\alpha }u|_{\beta }+C|u|_{0} \\
&=&\sup_{k}|\partial ^{\alpha }(\eta _{k}u)-E_{k,1}u-F_{k,1}u|_{\beta
}+C|u|_{0}.
\end{eqnarray*}

By (\ref{f20}) and (\ref{f2}), for each $\varepsilon >0$ there is a constant 
$C=C(\varepsilon ,\alpha ,\beta ,d,\delta )$ such that for every $u\in
C^{\alpha +\beta }(H)$%
\begin{equation*}
|u|_{\alpha ,\beta }\leq \sup_{k}|\partial ^{\alpha }\left( u\eta
_{k}\right) |_{\beta }+\varepsilon |\partial ^{\alpha }u|_{\beta }+C|u|_{0}.
\end{equation*}%
Therefore,%
\begin{equation*}
|u|_{\alpha ,\beta }\leq 2\sup_{k}|u\eta _{k}|_{\alpha ,\beta }+C|u|_{0},
\end{equation*}%
where the constant $C=C(\alpha ,\beta ,d,\delta )$. This estimate, together
with Lemma \ref{rem1}, implies%
\begin{equation}
|u|_{\alpha +\beta }\leq C_{1}\sup_{k}|u\eta _{k}|_{\alpha +\beta
}+C_{2}|u|_{0},  \label{f4}
\end{equation}%
where the constants $C_{1}=C_{1}(\alpha ,\beta ,d),C_{2}=C_{2}(\alpha ,\beta
,d,\delta ).$

Let $u\in C^{\alpha +\beta }(H)$ be a solution of (\ref{eq1}). Then $\eta
_{k}u$ satisfies the equation%
\begin{eqnarray}
\partial _{t}(\eta _{k}u) &=&A_{k}(\eta _{k}u)-\lambda (\eta _{k}u)+\eta
_{k}(Au-A_{k}u)+\eta _{k}f+\eta _{k}Bu  \label{eqk} \\
&&-F_{k}u-E_{k}u.  \notag
\end{eqnarray}%
By Theorem \ref{t1},%
\begin{equation*}
|\eta _{k}u|_{\alpha +\beta }\leq C\left( |\eta _{k}(Au-A_{k}u)|_{\beta
}+|\eta _{k}Bu|_{\beta }+|\eta _{k}f|_{\beta }+|F_{k}u|_{\beta
}+|E_{k}u|_{\beta }\right) ,
\end{equation*}%
where the constant $C=C(\alpha ,\beta ,d,\mu ,K,T)$. By (\ref{f4}),%
\begin{equation}
|u|_{\alpha +\beta }\leq C_{1}\left( \sup_{k}|\eta _{k}f|_{\beta }+I\right)
+C_{2}|u|_{0},  \label{form00}
\end{equation}%
where the constants $C_{1}=C_{1}(\alpha ,\beta ,d,\mu
,K,T),C_{2}=C_{2}(\alpha ,\beta ,d,K,\delta ,T)$ and%
\begin{equation*}
I=\sup_{k}\left( |\eta _{k}(Au-A_{k}u)|_{\beta }+|\eta _{k}Bu|_{\beta
}+|F_{k}u|_{\beta }+|E_{k}u|_{\beta }\right) .
\end{equation*}%
By Corollary \ref{cor1}, there is $\beta ^{\prime }\in (0,\beta )$ such that 
\begin{equation*}
|\eta _{k}(Au-A_{k}u)|_{\beta }\leq C_{1}\left( \delta ^{\beta }|u|_{\alpha
+\beta }+C_{2}|u|_{\alpha +\beta ^{\prime }}\right) ,
\end{equation*}%
where the constants $C_{1}=C_{1}(\alpha ,\beta ,d,K),C_{2}=C_{2}(\alpha
,\beta ,d,K,\delta )$. Therefore, for each $\varepsilon >0$ we can choose $%
\delta >0$ so that%
\begin{equation*}
|\eta _{k}(Au-A_{k}u)|_{\beta }\leq \varepsilon |u|_{\alpha +\beta
}+C|u|_{0},
\end{equation*}%
where the constant $C=C(\alpha ,\beta ,d,K,\varepsilon )$. Hence, by (\ref%
{f20}), (\ref{f2}) and Lemma \ref{l7}, for each $\varepsilon >0\,$\ we can
choose $\delta >0$ such that 
\begin{equation}
I\leq \varepsilon |u|_{\alpha +\beta }+C|u|_{0},  \label{f9}
\end{equation}%
where the constant $C=C(\alpha ,\beta ,d,K,\varepsilon )$. This estimate,
together with (\ref{form00}) and (\ref{f0}), implies  
\begin{equation}
|u|_{\alpha +\beta }\leq C[|f|_{\beta }+|u|_{0}],  \label{form01}
\end{equation}%
where the constant $C=C(\alpha ,\beta ,d,K,\mu ,T)$.

\bigskip On the other hand, according to (\ref{eqk}) and Theorem \ref{t1},%
\begin{eqnarray*}
|u|_{0} &\leq &\sup_{k}|\eta _{k}u|_{\beta }\leq \mu (\lambda
)\sup_{k}[|f|_{\beta }+|\eta _{k}(Au-A_{k})|_{\beta }+|\eta _{k}Bu|_{\beta }
\\
&&+|F_{k}u|_{\beta }+|E_{k}u|_{\beta }],
\end{eqnarray*}%
where $\mu (\lambda )\rightarrow 0$ as $\lambda \rightarrow \infty $. So,  
\begin{equation}
|u|_{0}\leq C\mu (\lambda )(|f|_{\beta }+|u|_{\alpha +\beta }).
\label{form02}
\end{equation}%
The inequalities (\ref{form01}) and (\ref{form02}) imply that there is $%
\lambda _{0}>0$ and a constant $C$ not depending on $u$ such that%
\begin{equation}
|u|_{\alpha +\beta }\leq C|f|_{\beta }  \label{form03}
\end{equation}%
if $\lambda \geq \lambda _{0}.$ If  $u\in C^{\alpha +\beta }(H)$ solves (\ref%
{eq1}) with $\lambda \leq \lambda _{0}$, then $\tilde{u}(t,x)=e^{-\lambda
(\lambda _{0}-\lambda )t}u(t,x)$ solves the same equation with $\lambda
=\lambda _{0}$, and by (\ref{form03})%
\begin{equation*}
|u|_{\alpha +\beta }\leq e^{(\lambda _{0}-\lambda )T}|\tilde{u}|_{\alpha
+\beta }\leq Ce^{(\lambda _{0}-\lambda )T}|f|_{\beta }.
\end{equation*}%
So, (\ref{form03}) holds for all $\lambda \geq 0$. 

By Theorem \ref{t1} and (\ref{f4}), there is a constant $C$ such that for
all $s\leq t\leq T,$ 
\begin{eqnarray*}
|u(t,\cdot )-u(s,\cdot )|_{\alpha /2+\beta } &\leq &C\sup_{k}|\eta
_{k}u(t,\cdot )-\eta _{k}u(s,\cdot )|_{\alpha /2+\beta } \\
&\leq &C(t-s)^{1/2}\left( |f|_{\beta }+|u|_{\alpha +\beta }\right) .
\end{eqnarray*}%
Therefore there is a constant $C$ such that for all $s\leq t\leq T,$%
\begin{equation*}
|u(t,\cdot )-u(s,\cdot )|_{\alpha /2+\beta }\leq C(t-s)^{1/2}|f|_{\beta }.
\end{equation*}

We finish the proof applying the continuation by parameter argument. Let 
\begin{equation*}
\,L_{\tau }u=\tau Lu+\left( 1-\tau \right) \partial ^{\alpha }u,\tau \in 
\left[ 0,1\right] .
\end{equation*}
We introduce the space $\hat{C}^{\alpha +\beta }\left( H\right) $ of
functions $u\in C^{\alpha +\beta }(H)$ such that for each $\left( t,x\right) 
$, $u\left( t,x\right) =\int_{0}^{t}F\left( s,x\right) \,ds,$where $F\in
C^{\beta }\left( H\right) .$ It is a Banach space with respect to the norm 
\begin{equation*}
|\left\vert u\right\vert |_{\alpha ,\beta }=\left\vert u\right\vert _{\alpha
+\beta }+\left\vert F\right\vert _{\beta }.
\end{equation*}%
Consider the mappings $T_{\tau }:\hat{C}^{\alpha +\beta }\left( H\right)
\rightarrow C^{\beta }(H)$ defined by%
\begin{equation*}
u\left( t,x\right) =\int_{0}^{t}F\left( s,x\right) \,ds\longmapsto F-L_{\tau
}u.
\end{equation*}%
Obviously, for some constant $C$ not depending on $\tau ,$ 
\begin{equation*}
\left\vert T_{\tau }u\right\vert _{\beta }\leq C|\left\vert u\right\vert
|_{\alpha ,\beta }.
\end{equation*}%
On the other hand, there is a constant $C$ not depending on $\tau $ such
that for all $u\in \hat{C}^{\alpha +\beta }\left( H\right) $%
\begin{equation}
|\left\vert u\right\vert |_{\alpha ,\beta }\leq C\left\vert T_{\tau
}u\right\vert _{\beta }.  \label{cp1}
\end{equation}%
Indeed, 
\begin{equation*}
u\left( t,x\right) =\int_{0}^{t}F\left( s,x\right) \,ds=\int_{0}^{t}\left(
L_{\tau }u+(F-L_{\tau }u)\right) (s,x)\,ds,
\end{equation*}%
and, according to the estimate (\ref{form03}), there is a constant $C$ not
depending on $\tau $ such that 
\begin{equation}
\left\vert u\right\vert _{\alpha +\beta }\leq C\left\vert T_{\tau
}u\right\vert _{\beta }=C\left\vert F-L_{\tau }u\right\vert _{\beta }.
\label{cp2}
\end{equation}%
Thus,%
\begin{eqnarray*}
|\left\vert u|\right\vert _{\alpha ,\beta } &=&\left\vert u\right\vert
_{\alpha +\beta }+\left\vert F\right\vert _{\beta }\leq \left\vert
u\right\vert _{\alpha +\beta }+\left\vert F-L_{\tau }u\right\vert _{\beta
}+\left\vert L_{\tau }u\right\vert _{\beta } \\
&\leq &C\left( \left\vert u\right\vert _{\alpha +\beta }+\left\vert
F-L_{\tau }u\right\vert _{\beta }\right) \leq C\left\vert F-L_{\tau
}u\right\vert _{\beta }=C\left\vert T_{\tau }u\right\vert _{\beta },
\end{eqnarray*}%
and (\ref{cp1}) follows. Since $T_{0}$ is an onto map, by Theorem 5.2 in 
\cite{GiT83} all the $T_{\tau }$ are onto maps and the statement follows.

\section{Martingale problem}

In this section, we consider the martingale problem associated with the
operator%
\begin{equation*}
L^{0}=A+B^{0}\text{,}
\end{equation*}%
where $B^{0}\,$\ is the operator $B$ defined by (\ref{for2}) with $\rho \geq
0$ and $l=0.$

Let $D=D([0,T],\mathbf{R}^{d})$ be the Skorokhod space of cadlag $\mathbf{R}%
^{d}$-valued trajectories and let $X_{t}=X_{t}(w)=w_{t},w\in D,$ be the
canonical process on it. 

Let%
\begin{equation*}
\mathcal{D}_{t}=\sigma (X_{s},s\leq t),\mathcal{D}=\vee _{t}\mathcal{D}_{t},%
\mathbb{D=}\left( \mathcal{D}_{t+}\right) ,t\in \lbrack 0,T].
\end{equation*}%
We say that a probability measure $\mathbf{P}$ on $\left( D,\mathcal{D}%
\right) $ is a solution to the $(s,x,L)$-martingale problem (see \cite{st}, 
\cite{MiP923}) if $\mathbf{P}(X_{r}=x,0\leq r\leq s)=1$ and for all $u\in
C_{0}^{\infty }(H)$ the process%
\begin{equation}
u(t,X_{t})-\int_{s}^{t}[\partial _{t}u(r,X_{r})+L^{0}u(r,X_{r})]dr
\label{mart1}
\end{equation}
is a $(\mathbb{D},\mathbf{P})$-martingale. We denote $\mathcal{S}(s,x,L^{0})$
the set of all solutions to the problem $(s,x,L^{0})$-martingale problem.

\begin{lemma}
\label{ml1}let $\alpha \in (0,2),\beta \in (0,1]$ and Assumptions A, B1 and
B2 with $\rho \geq 0$ and $l=0$ be satisfied. Let $\mathbf{P\in }\mathcal{S}%
(s,x,L^{0}),f\in C^{\beta }(H),$ and let $u\in C^{\alpha +\beta }(H)$ be a
solution to the Cauchy problem%
\begin{eqnarray}
\partial _{t}u(t,x)+L_{t}^{0}u(t,x) &=&f(t,x),(t,x)\in H,  \label{mart0} \\
u(T,x) &=&0.  \notag
\end{eqnarray}%
Then the process (\ref{mart1}) is a ($\mathbb{D},\mathbf{P})$-martingale and%
\begin{equation}
u(s,x)=-\mathbf{P}_{s,x}\int_{s}^{T}f(r,X_{r})dr,(s,x)\in H.  \label{mart2}
\end{equation}
\end{lemma}

\begin{proof}
Let $\zeta _{\varepsilon }$ be the function introduced in the proof of
Theorem \ref{t1} and 
\begin{equation*}
u_{\varepsilon }(t,x)=u(t,\cdot )\ast \zeta _{\varepsilon }(x),(t,x)\in
H,\varepsilon \in (0,1).
\end{equation*}
Let $r<t$ and $h$ be a bounded $\mathcal{F}_{r}$-measurable random variable$.
$Then%
\begin{equation*}
\mathbf{P}_{s,x}\{h[u_{\varepsilon }(t,X_{t})-u_{\varepsilon
}(r,X_{r})-\int_{r}^{t}(\partial _{t}u_{\varepsilon
}(s,X_{s})+L^{0}u_{\varepsilon }(s,X_{s}))ds]\}=0.
\end{equation*}%
Passing to the limit in this equality as $\varepsilon \rightarrow 0$ and
using (\ref{mart0}), we get%
\begin{equation*}
\mathbf{P}_{s,x}\{h[u(t,X_{t})-u(r,X_{r})-\int_{r}^{t}f(s,X_{s}))ds]\}=0.
\end{equation*}%
In particular, for $t=T,r=0,h=1,$ (\ref{mart2}) follows.
\end{proof}

\begin{proposition}
\label{prop3}Let Assumptions A, B1 and B2 with $\rho \geq 0$ and $l=0$ be
satisfied. Then for each $(s,x)\in H$ there is a unique solution $\mathbf{P}%
_{s,x}$ to the martingale problem $(s,x,L^{0}),$ and the process $\left(
X_{t},\mathbb{D},(\mathbf{P}_{s,x})\right) $ is strong Markov.

If, in addition,%
\begin{equation*}
\lim_{R\rightarrow \infty }\int_{0}^{T}\sup_{x}\int_{|c(t,x,\upsilon
)|>R}\rho (t,x,\upsilon )\pi (d\upsilon )dt=0,
\end{equation*}%
then the function $\mathbf{P}_{s,x}$ is weakly continuous in $(s,x)$.
\end{proposition}

\begin{proof}
Since the coefficients of $L^{0}$ are H\"{o}lder continuous, it follows by
Theorem IX.2.31 in \cite{jaks} that  the set $\mathcal{S}(s,x,L^{0})\neq
\emptyset .$ For $f\in C^{\beta }(H)$, let $u\,\in C_{\alpha +\beta }(H)$ be
the solution to (\ref{eq1}). By Lemma \ref{ml1},%
\begin{equation*}
u(s,x)=\mathbf{P}_{s,x}\int_{s}^{T}f(r,X_{r})dr,\mathbf{P}_{s,x}\in \mathcal{%
S}(s,x,L).
\end{equation*}%
Therefore, by Lemma 2.4 \cite{MiP923}, the measure $\mathbf{P}_{s,x}\in 
\mathcal{S}(s,x,L^{0})$ is unique. By Lemma 2.2 in \cite{MiP923}, the
process $\left( X_{t},\mathbb{D},(\mathbf{P}_{s,x})\right) $ is strong
Markov. The continuity of the function $(s,x)\rightarrow \mathbf{P}_{s,x}$
follows from Theorems IX.2.22 and IX.3.9 in \cite{jaks}.
\end{proof}

\end{document}